\documentclass[11pt]{article}
\usepackage{latexsym}
\usepackage{amsfonts}
\usepackage{amsmath}

\newtheorem{theorem}{Theorem}[section]
\newtheorem{conjecture}{Conjecture}[section]

\title{The Erd\"{o}s-Hajnal Conjecture---A Survey}
\author{Maria Chudnovsky
\thanks{Partially supported by NSF grants DMS-1001091 and IIS-1117631. }\\
Columbia University, New York NY 10027
}

\begin{document}
\maketitle
\begin{abstract}
The Erd\"{o}s-Hajnal conjecture states that for every graph $H$, there exists 
a constant $\delta(H) > 0$ such that every graph $G$ with no induced subgraph 
isomorphic to $H$  has either a clique or a stable set of size at least
$|V(G)|^{\delta(H)}$. This paper is a survey of some of the known results on this
conjecture.
\end{abstract}
\section{Introduction}

All graphs in this paper are finite and simple.  Let $G$ be a graph.
We denote by $V(G)$ the vertex set of $G$.
The {\em complement} $G^c$ of $G$ is the graph with vertex set $V(G)$,
such that two vertices are adjacent in $G$ if and only if they are non-adjacent
in $G^c$.  A {\em clique} in $G$ is a set of vertices all pairwise adjacent. 
A {\em stable set} in $G$ is a set of vertices all pairwise non-adjacent
(thus a stable set in $G$ is a clique in $G^c$). Given a graph $H$, we say
that $G$ is {\em $H$-free} if $G$ has no induced subgraph isomorphic to $H$.
For a family $\mathcal{F}$ of graphs, we say that $G$ is $\mathcal{F}$-free
if $G$ is $F$-free for every $F \in \mathcal{F}$.

It is a well-known theorem of  Erd\"{o}s \cite{Erdos1} that there exist graphs 
on $n$ vertices, with no clique  or stable set of size larger than $O(\log n)$. 
However, in 1989 Erd\"{o}s and  Hajnal \cite{EH} made a conjecture suggesting 
that the situation  is  dramatically different for graphs that are $H$-free 
for some fixed graph $H$, the following:
\begin{conjecture}
\label{EHconj}
For every graph $H$, there exists a constant $\delta(H) > 0$  such that every
$H$-free graph $G$ has either a clique or a stable set of size at least
$|V(G)|^{\delta(H)}$.
\end{conjecture}

This is the {\em Erd\"{o}s-Hajnal conjecture}. The same paper \cite{EH} also
contains  a partial result
toward Conjecture~\ref{EHconj}, showing that $H$-free graphs behave very differently from 
general graphs:

\begin{theorem}
\label{EHthm}
For every graph $H$, there exists a constant $c(H) > 0$ such that every
$H$-free graph $G$ has either a clique or a stable set of size at least
$e^{c(H) \sqrt{log |V(G)|}}$.
\end{theorem}

However, obtaining the polynomial bound of Conjecture~\ref{EHconj} seems to be a lot 
harder, and Conjecture~\ref{EHconj} is still open. The goal of this paper is to survey 
some recent results on Conjecture~\ref{EHconj}.

We start with some definitions.
We say that a graph $H$ has the 
{\em Erd\"{o}s-Hajnal property} if there 
exists a constant $\delta(H) > 0$ such that every
$H$-free graph $G$ has either a clique or a stable set of size at least
$|V(G)|^{\delta(H)}$. Clearly, $H$ has the Erd\"{o}s-Hajnal property if and
only if $H^c$ does.

Let $G$ be a graph. For $X \subseteq V(G)$, we denote by $G|X$ the subgraph of 
$G$ induced by $X$. We write $G \setminus X$ for $G|(V(G) \setminus X)$,
and $G \setminus v$ for $G \setminus \{v\}$, where $v \in V(G)$.
We  denote by $\omega(G)$ the maximum size of a clique in 
$G$,  by $\alpha(G)$ the maximum size of a stable set in $G$, and by $\chi(G)$ 
the chromatic number of $G$. The graph $G$ is   {\em perfect} if  
$\chi(H)=\omega(H)$ for every induced subgraph $H$ of $G$. The Strong Perfect
Graph Theorem \cite{CRST} characterizes perfect graphs by forbidden induced 
subgraphs:

\begin{theorem}
\label{SPGT}
A graph $G$ is perfect if and only if no induced subgraph of $G$ or $G^c$ is 
an odd cycle of length at least five.
\end{theorem}

Thus, in a perfect graph $G$
$$|V(G)| \leq \chi(G) \alpha(G) = \omega(G) \alpha(G),$$
and so 

\begin{theorem}
\label{perfect}
If $G$ is perfect, then either $\omega(G) \geq \sqrt{|V(G)|}$, or 
$\alpha(G)  \geq \sqrt{|V(G)|}$.
\end{theorem}

This turns out to be a useful observation in the study of Conjecture~\ref{EHconj}; in
fact it is often convenient to work with the following equivalent version of
Conjecture~\ref{EHconj}:

\begin{conjecture}
\label{EHconj3}
For every graph $H$, there exists a constant $\psi(H) > 0$, such that every
$H$-free graph $G$ has a perfect induced subgraph  with at  least
$|V(G)|^{\psi(H)}$ vertices.
\end{conjecture}

The equivalent of Conjecture~\ref{EHconj} and Conjecture~\ref{EHconj3} follows from Theorem~\ref{perfect}.
The main advantage of Conjecture~\ref{EHconj3} is that instead of having two outcomes: 
a large clique or a large stable set, it only has one, namely a large perfect
induced subgraph, thus making inductive proofs easier.

This paper is organized as follows. In Section~2 we discuss graphs and
families of graphs that are known to have the Erd\"{o}s-Hajnal property.
Section~3 deals with weakenings of  Conjecture~\ref{EHconj} that are known to be true
for all graphs. Is Section~4 we state an analogue of Conjecture~\ref{EHconj} for 
tournaments, and discuss related  results and techniques. Sections~5 and 6
deal with special cases of Conjecture~\ref{EHconj} and its tournament analogue when
we restrict our attention to graphs $H$ with a certain value of $\delta(H)$.

\section{Graphs with the Erd\"{o}s-Hajnal property} 

Obviously, Conjecture~\ref{EHconj} can be restated as follows:

\begin{conjecture}
\label{EHconj2}
Every graph has the Erd\"{o}s-Hajnal property.
\end{conjecture}

However, at the moment, only very few graphs have been shown to have 
the  Erd\"{o}s-Hajnal property. The goal of this section is to describe all 
such graphs. 
It is clear that graphs on at most two vertices have the property. 
Complete graphs and their complements have the property; this follows from the 
famous Ramsey theorem. If $H$ is the two-edge path, then every
$H$-free graph $G$ is the disjoint union of cliques, and thus has either
a clique or a stable set of size $\sqrt{|V(G)|}$, so the two-edge path has
the property. By taking complements, this shows that all three-vertex 
graphs have the property.

Another graph for which the Erd\"{o}s-Hajnal property is easily established
is the three-edge path. It follows immediately from Theorem~\ref{SPGT} that all
graphs with no induced subgraph isomorphic to the three-edge path are perfect;
this fact  can also be obtained by an easy induction from the following 
theorem of  Seinsche \cite{Seinsche}:

\begin{theorem}
\label{3edgepath}
If $G$ is a graph with at least two vertices, and  no induced subgraph of $G$ 
is isomorphic to the three-edge path, then either $G$ or $G^c$ is not connected.
\end{theorem}

Let us now define the {\em substitution} operation. Given graphs $H_1$ and $H_2$,
on disjoint vertex sets, each with at least two vertices, and $v \in V(H_1)$, 
we say that $H$ is {\em obtained 
from $H_1$ by substituting $H_2$ for $v$}, or {\em obtained from $H_1$ and $H_2$
by substitution} (when the details are not important)  if:
\begin{itemize}
\item $V(H)=(V(H_1) \cup V(H_2)) \setminus \{v\}$,
\item $H|V(H_2)=H_2$,
\item $H|(V(H_1) \setminus \{v\})=H_1 \setminus v$, and
\item $u \in V(H_1)$ is adjacent in $H$ to $w \in V(H_2)$ if and only if
$u$ is adjacent in $H_1$  to $v$.
\end{itemize}

A graph is {\em prime} if it is not obtained from smaller graphs by 
substitution.

In \cite{APS} Alon, Pach, and Solymosi proved that the Erd\"{o}s-Hajnal property
is preserved under substitution:
\begin{theorem}
\label{substitution}
If $H_1$ and $H_2$ are graphs with the Erd\"{o}s-Hajnal property, and $H$
is obtained from $H_1$ and $H_2$ by substitution, then $H$ has  the 
Erd\"{o}s-Hajnal property.
\end{theorem}

This is the only operation known today that allows us to build bigger graphs
with the Erd\"{o}s-Hajnal property from smaller ones.
The idea of the proof of Theorem~\ref{substitution} is to notice the following: 
since both $H_1$ and $H_2$ have the  Erd\"{o}s-Hajnal property, it follows that 
if an $H$-free graph $G$ does not contain a ``large'' clique or stable set, 
then  every induced subgraph of $G$ with at least $|V(G)|^{\epsilon}$
vertices (where $\epsilon$ depends on the precise definition of ``large'') 
contains an induced copy of $H_1$ and an induced copy of $H_2$.  Let 
$v \in V(H_1)$ be such that $H$ is obtained from $H_1$ by substituting $H_2$ 
for $v$. Then counting shows that some copy of $H_1\setminus v$ in $G$ can be 
extended to $H_1$ in at least $n^{\epsilon}$ ways. But this guarantees that there
is a copy of $H_2$ among the possible extensions, contrary to the fact that $G$
is $H$-free, and Theorem~\ref{substitution} follows.

Since the only prime graph on four vertices is the three-edge path,
using Theorem~\ref{substitution} and the fact that the three-edge-path has
the Erd\"{o}s-Hajnal property, it is easy to check that all graphs on at most
four vertices have the Erd\"{o}s-Hajnal property. Moreover, there are only
four prime graphs on five vertices:
\begin{itemize}
\item the cycle of length five
\item the four-edge path
\item the complement of the four-edge path
\item the bull (the graph with vertex set $\{a_1,a_2,a_3,b_1,b_2\}$ and
edge set $\{a_1a_2,a_2a_3,a_1a_3, a_1b_1, a_2b_2\}$),
\end {itemize}
and Theorem~\ref{substitution} implies that all the other graphs on  five
vertices have the Erd\"{o}s-Hajnal property.

A much harder argument \cite{CSafra} shows that
\begin{theorem}
\label{EHbull}
The bull has the  the Erd\"{o}s-Hajnal property. Moreover, every
bull-free graph $G$ has a clique or a stable set of size at least
$|V(G)|^{\frac{1}{4}}$.
\end{theorem}

The exponent $\frac{1}{4}$ is in fact best possible because of the following 
construction. Take a triangle-free graph $T$ with $m$ vertices and no
stable set of size larger than $\sqrt{m \log m}$ (such graphs exist by an old 
theorem of Kim \cite{Kim}). Let $G$ be obtained from $T$ by
substituting a copy of $T^c$ for every vertex of $T$. Then $|V(G)|=m^2$,
and it  is easy to check that  $G$ has no clique or stable set of size
larger than $2\sqrt{m \log m}$.

The proof of Theorem~\ref{EHbull} uses  structural methods to show that 
every prime bull-free graph belongs to a certain
subclass where a clique or a stable set of the appropriate size can be shown 
to exist. For graphs that are not prime, the result follows by induction.

Let us say that a function $f: V(G) \rightarrow [0,1]$ is {\em good} if
for every perfect induced subgraph $P$ of $G$
$$\Sigma_{v \in V(P)}f(v) \leq 1.$$ 
For $\alpha \geq 1$, the graph $G$ is {\em $\alpha$-narrow} if for every 
good function $f$
$$\Sigma_{v \in V(G)}f(v)^{\alpha} \leq 1.$$ 

Thus perfect graphs are $1$-narrow.
Let $G$ be an  $\alpha$-narrow graph for some  $\alpha \geq 1$, and let 
$K=max |V(P)|$ where the maximum is taken over all perfect induced subgraphs of
$G$. Then the function $f(v)=\frac{1}{K}$ (for all $v \in V(G)$) is good,
and so, since $G$ is $\alpha$-narrow,  $\frac{|V(G)|}{K^{\alpha}} \leq 1$.
Thus $K \geq |V(G)|^{\frac{1}{\alpha}}$. By Theorem~\ref{perfect}, this implies 
that in 
order to prove that a certain graph $H$ has the Erd\"{o}s-Hajnal property, it
is enough to show that there exists $\alpha \geq 1$ such that all $H$-free 
graphs are $\alpha$-narrow. This conjecture was formally stated in 
\cite{CZwols}:
\begin{conjecture}
\label{narrow}
For every graph $H$, there exists a constant $\alpha(H) \geq 1$ such that every
$H$-free graph $G$ is $\alpha(H)$-narrow.
\end{conjecture}

Fox \cite{Fox} (see \cite{CSeymour2} for details)  proved that 
Conjecture~\ref{narrow} is in fact 
equivalent to Conjecture~\ref{EHconj}.
More specifically, he showed that for every graph $H$ with  the 
Erd\"{o}s-Hajnal property, there
exists a constant $\alpha(H)$ such that  every $H$-free graph is 
$\alpha(H)$-narrow.
However, at least for the purely structural approach, 
Conjecture~\ref{narrow} seems to be more convenient to work with than 
\ref{EHconj}. In fact, what is really
proved in \cite{CSafra} is the following:
\begin{theorem}
\label{narrowbull}
Every bull-free graph is $2$-narrow.
\end{theorem}

And the inductive step proving Theorem~\ref{narrowbull} for bull-free graphs that are
not prime is that $\alpha$-narrowness is preserved
under substitutions:

\begin{theorem}
\label{subnarrow}
If $H_1$ and $H_2$ are $\alpha$-narrow for $\alpha \geq 1$, and
$H$ is obtained from $H_1$ and $H_2$ by substitution, then $H$ is also
$\alpha$-narrow. 
\end{theorem}

Conjecture~\ref{EHconj} is still open for the four-edge path, its complement, and the
five-cycle;  and  no prime graph on at least six vertices is known to have
the Erd\"{o}s-Hajnal property.

We remark that the bull is a self-complementary graph, and
one might think that to be the reason for its better behavior.
This philosophy is supported by the following result:
\begin{theorem}
\label{P4P4C}
Every graph with no induced subgraph isomorphic to the four-edge-path
or the complement of the four-edge-path is $2$-narrow.
\end{theorem}
This follows from Theorem~\ref{subnarrow} and  from 
(a restatement of) a theorem of Fouquet \cite{Fouquet}:
\begin{theorem}
\label{Fouquet}
Every prime graph with no induced subgraph isomorphic to the four-edge-path
or the complement of the  four-edge-path is either perfect or isomorphic to 
the five-cycle.
\end{theorem}

On the other hand,  the five-cycle is another self-complementary graph,
and yet it seems to be completely intractable. We thus propose the following
conjecture that may be slightly easier than the full Conjecture~\ref{EHconj} in
special cases:
\begin{conjecture}
\label{EHconj4}
For every graph $H$, there exists a constant $\epsilon(H) > 0$, such that every
$\{H,H^c\}$-free graph $G$ has either a clique or a stable set of size at least
$|V(G)|^{\epsilon(H)}$.
\end{conjecture}

In \cite{CSeymour2} Conjecture~\ref{EHconj4} was proved in the case when $H$ is the
five-edge path. This is a strengthening of a result of \cite{CZwols}.

\section{Approximate results}

The previous section listed a few graphs for which Conjecture~\ref{EHconj2} is known to 
hold. The goal of this section is to list facts that are true for all
graphs, but that do not achieve the full strength conjectured in Conjecture~\ref{EHconj2}.
The first such statement is Theorem~\ref{EHthm} which we already mentioned in the 
Introduction. 

For two disjoint subsets $A$ and $B$ of $V(G)$, we say that $A$ is 
{\em complete} to $B$ if every vertex of $A$ is adjacent to every vertex of
$B$, and that  $A$ is  {\em anticomplete} to $B$ if every vertex of $A$ is
non-adjacent to every vertex of $B$. If $A=\{a\}$ for some $a \in V(G)$,
we write ``$a$ is complete (anticomplete) to $B$'' instead of  ``$\{a\}$ is
complete (anticomplete) to $B$''. Here is another theorem similar 
to Theorem~\ref{EHthm}, due to Erd\"{o}s, Hajnal and Pach \cite{EHP}.

\begin{theorem}
\label{twosets}
For every graph $H$, there exists a constant $\delta(H)>0$ such that
for every $H$-free graph $G$ there exist two disjoint subsets 
$A,B \subseteq V(G)$ with the following properties:
\begin{enumerate}
\item $|A|,|B| \geq |V(G)|^{\delta(H)}$, and
\item either $A$ is complete to $B$, or $A$ is anticomplete to $B$.
\end{enumerate}
\end{theorem}
 
The idea of the proof here is to partition $V(G)$ into $|V(H)|$ equal 
subsets (which we call ``sets of candidates''), and then try to build an 
induced copy of $H$ in $G$, one 
vertex at a time, where each  vertex of $H$ is chosen from the corresponding
set of candidates.  In this process, sets of candidates shrink 
at every step, and since $G$ is $H$-free, eventually we reach a situation 
where there do not exist enough vertices in one of the sets of candidates with 
the right  adjacencies to another. At this stage we obtain the sets $A$ and 
$B$, as required in Theorem~\ref{twosets}.

Theorem~\ref{twosets} was recently strengthened by Fox and Sudakov in \cite{FS}:
\begin{theorem}
\label{halfway}
For every graph $H$, there exists a constant $\delta(H)>0$ such that
for every $H$-free graph $G$ with $\omega(G) < {{|V(G)|}^{\delta(H)}}$ 
there exist two disjoint subsets 
$A,B \subseteq V(G)$, with the following properties:
\begin{enumerate}
\item $|A|,|B| \geq |V(G)|^{\delta(H)}$, and
\item $A$ is anticomplete to $B$.
\end{enumerate}
\end{theorem}

In \cite{LRSTT} another weakening of Conjecture~\ref{EHconj2} is considered. It is shown 
that for every $H$, the proportion of $H$-free graphs with 
$n$ vertices and no ``large'' cliques or stable sets tends to zero as 
$n \rightarrow \infty$. Let $\mathcal{F}_H^n$ be the class
of all $H$-free graphs on $n$ vertices. Let  $\mathcal{Q}_H^{n,\epsilon}$
be the subclass of  $\mathcal{F}_H^n$, consisting of all graphs $G$
that have either a clique or a stable set of size at least
$n^{\epsilon}$.
The main result of \cite{LRSTT} is the following:
\begin{theorem}
\label{almostall}
For every graph $H$, there exists a constant $\epsilon(H)>0$ such that 
${\frac {|\mathcal{Q}_H^{n,\epsilon(H)}|}{|\mathcal{F}_H^n|}} \rightarrow 1$
as $n \rightarrow \infty$.
\end{theorem}
The proof of Theorem~\ref{almostall} involves an application of Szemer\'{e}di's 
Regularity Lemma \cite{Szemeredi}. Also, somewhat surprisingly, it uses 
Theorem~\ref{EHbull}.

\section{Tournaments}

A {\em tournament} is a directed graph, where for every two vertices $u,v$ 
exactly one of the (ordered) pairs $uv$ and $vu$ is an edge. A tournament is
{\em transitive} if it has no directed cycles (or, equivalently, no cyclic 
triangles). For a tournament $T$, we denote by $\alpha(T)$ the maximum number of
vertices of a transitive subtournament of $T$. Transitive subtournaments seem 
to be a good analogue of both cliques and stable sets in graphs; furthermore, 
like  induced perfect subgraphs in Conjecture~\ref{EHconj3}, transitive 
tournaments have the advantage of being one object instead of two.

For tournaments $S$ and $T$, we say that $T$ is {\em $S$-free} if no
subtournament of $T$ is isomorphic to $S$. As with graphs, if $\mathcal{S}$
is a family of tournaments, then $T$ is {\em $\mathcal{S}$-free} if $T$ is 
$S$-free for every $S \in \mathcal{S}$. In \cite{APS} the following conjecture
was formulated, and shown to be equivalent to Conjecture~\ref{EHconj}:
\begin{conjecture}
\label{EHtourn}
For every tournament $S$, there exists a constant $\delta(S) > 0$ such that 
every $S$-free tournament  $T$ satisfies $\alpha(T) \geq |V(T)|^{\delta(H)}$.
\end{conjecture}
As with graphs, let us say that a tournament $S$  has the 
{\em Erd\"{o}s-Hajnal property} if there exists $\delta(S) > 0$ such that 
every $S$-free tournament  $T$ satisfies $\alpha(T) \geq |V(T)|^{\delta(H)}$.
We remark that, like in graphs, the maximum number of vertices of a  
transitive subtournament in a  random $n$-vertex tournament is $O(\log n)$ \cite{Erdos1}.

A {\em substitution} operation can be defined for tournaments as follows.
Given  tournaments $S_1$ and $S_2$, with disjoint vertex sets and each with 
at least two vertices,  and a vertex $v \in V(S_1)$, 
we say that $S$ is {\em obtained from $S_1$ by substituting $S_2$ for $v$}
(or just {\em obtained by substitution from $S_1$ and $S_2$}) if
$V(S)=(V(S_1) \cup V(S_2)) \setminus \{v\}$ and $uw$ is an edge of $S$ if and
only if one of the following holds:
\begin{itemize}
\item $u,w \in V(S_1)$ and $uw$ is an edge of $S_1$
\item $u,w \in V(S_2)$ and $uw$ is an edge of $S_2$
\item $u \in S_1, w \in S_2$ and $uv$ is an edge of $S_1$
\item $u \in S_2$, $w \in S_1$, and $vw$ is an edge of $S_1$.
\end{itemize}
A tournament is {\em prime} if it is not obtained by substitution from
smaller tournaments. Repeating the proof of Theorem~\ref{substitution} in the
setting of tournaments instead of graphs,  it is easy to 
show that
\begin{theorem}
\label{subsitutiontourn}
If $S_1$ and $S_2$ are tournaments with the Erd\"{o}s-Hajnal property, and $S$
is obtained from $S_1$ and $S_2$ by substitution, then $S$ has  the 
Erd\"{o}s-Hajnal property.
\end{theorem}

Clearly, all tournaments on at most three vertices have the  Erd\"{o}s-Hajnal 
property, and it is easy to check that there are no prime four-vertex 
tournaments. Consequently, all four-vertex tournaments have the 
Erd\"{o}s-Hajnal  property. So far this is very similar to the state of 
affairs in graphs, but here is a fact to which we do not have a graph analogue:
we can define an infinite family of prime tournaments all with the
Erd\"{o}s-Hajnal  property (recall that the largest prime graph known to
have the property is the bull). Let us describe this family.

First we need some definitions.
Let $T$ be a tournament, and let $(v_1, \ldots, v_{|T|})$ be an ordering of its 
vertices; denote it by $\theta$. We say that an edge $v_{j}v_{i}$ of $T$ is a 
{\em backward edge}  under $\theta$  if $i<j$.  
The {\em graph of backward edges}  under $\theta$,
denoted by $B(T, \theta)$, has vertex set $V(T)$,  and  
$v_i v_j \in E(B(T, \theta))$ if and only if  $v_iv_j$ or 
$v_jv_i$ is a backward edge of $T$ under the ordering $\theta$. 
For an integer $t>0$, we call the graph $K_{1,t}$ a {\em star}. Let $S$ be a
star with vertex set $\{c, l_1, \ldots, l_t\}$, where $c$ is adjacent to $l_1, \ldots, l_t$ and $\{l_1, \ldots, l_t\}$ is a stable set. We call $c$ the 
{\em center of the star}, and
$l_1, \ldots, l_t$ {\em the leaves of the star}. A {\em right star}
in $B(T, \theta)$ is an induced subgraph with vertex set 
$\{v_{i_0}, \ldots, v_{i_t}\}$, such that 
$B(T,\theta)|\{v_{i_0}, \ldots, v_{i_t}\}$ is a star with center $v_{i_t}$, 
and  $i_t > i_0, \ldots, i_{t-1}$.   A {\em left star}
in $B(T, \theta)$ is an induced subgraph with vertex set 
$\{v_{i_0}, \ldots, v_{i_t}\}$, such that 
$B(T,\theta)|\{v_{i_0}, \ldots, v_{i_t}\}$ is a star with center $v_{i_0}$, 
and  $i_0 < i_1, \ldots, i_t$.   
Finally, a {\em star}   in $B(T,\theta)$,  is a left star or a right star.
A tournament $T$ is a {\em galaxy} if there exists an ordering $\theta$ of its vertices such that every component of $B(T, \theta)$ is  a star or a 
singleton, and
\begin{itemize}
\item no center of a star appears in $\theta$  between two leaves 
of another star.
\end{itemize}

In \cite{BCC} the following is proved:
\begin{theorem}
\label{galaxy}
Every galaxy has the Erd\"{o}s-Hajnal property.
\end{theorem}

The proof uses the directed version of Szemer\'{e}di's Regularity Lemma 
formulated  in \cite{AS},  and extensions of ideas from the proof 
of Theorem~\ref{twosets}. Instead of  starting with arbitrary sets of candidates, the 
way it is done in  Theorem~\ref{twosets}, we get a head start  by using sets given by 
a regular partition.  Let us describe the proof  in a little more detail.

Let $T$ be a tournament.
For disjoint subsets $A,B$  of $V(T)$, we say that $A$ is {\em complete to} 
$B$ if every vertex of $A$ is adjacent to every vertex of $B$. We say that $A$ is {\em complete from} $B$ if $B$ is complete to $A$.
Denote by $e_{A,B}$ the number of directed edges $ab$, where $a \in A$ and 
$b \in B$. We define the {\em directed density from A to B} to be
$d(A,B)=\frac{e_{A,B}}{|A||B|}.$ 

Given $\epsilon>0$ we call a pair $(X,Y)$ of disjoint subsets of $V(T)$ 
$\epsilon$-{\em regular} if all $A \subseteq X$ and $B \subseteq Y$ with 
$|A| \geq \epsilon |X|$ and $|B| \geq \epsilon |Y|$ satisfy:
$|d(A,B)-d(X,Y)| \leq \epsilon$ and  $|d(B,A)-d(Y,X)| \leq \epsilon.$

Consider a partition $\{V_{0},V_{1},...,V_{k}\}$ of $V(T)$ in which one set $V_{0}$ has been singled out as an {\em exceptional} set. (This exceptional set $V_{0}$ may be empty). Such a partition is called an {\em $\epsilon$-regular partition of $T$} if it satisfies the following three conditions: 
\begin{itemize}
\item $|V_{0}| \leq \epsilon |V|$ 
\item $|V_{1}|=...=|V_{k}|$
\item all but at most $\epsilon k^{2}$ of the pairs $(V_{i},V_{j})$ with $1 \leq i < j \leq k$ are $\epsilon$-regular.
\end{itemize}

The following was proved in \cite{AS}:

\begin{theorem}
\label{regularitylemma}
For every $\epsilon>0$ and every $m \geq 1$ there exists an integer $DM=DM(m,\epsilon)$ such that every tournament of order at least $m$ admits an $\epsilon$-regular partition $\{V_{0},V_{1},...,V_{k}\}$ with $m \leq k \leq DM$.
\end{theorem}

The following is also a result from~\cite{AS}; it is the directed analogue of 
a well-known lemma for undirected graphs \cite{BT}.

\begin{theorem}
\label{universalgraphlemma}
Let $k \geq 1$ be an integer, and let $0 < \lambda < 1$. Then
there exists a constant $\eta_0$ (depending on $k$ and $\lambda$) with the
following property. Let  H be a tournament with vertex set 
$\{x_{1},...,x_{k}\}$, and  let $T$ be a tournament with vertex set $V(T)=\bigcup_{i=1}^{k} V_{i}$, where the $V_{i}$'s are disjoint sets, each of order at least one. Suppose that each pair $(V_{i},V_{j})$, $1 \leq i < j \leq k$ is $\eta$-regular, that $d(V_{i},V_{j}) \geq \lambda$ and $d(V_{j},V_{i}) \geq \lambda$. Then there exist vertices $v_{i} \in V_{i}$ for $i \in \{1, \ldots, k\}$, such that the map $x_{i} \rightarrow v_{i}$ gives an isomorphism between $H$ and the subtournament of $T$ induced by $\{v_{1},...,v_{k}\}$ provided that $\eta \leq \eta_{0}$.
\end{theorem}

Now, given a galaxy $G$, we start with a  regular partition of a $G$-free
tournament given by  Theorem~\ref{regularitylemma}. Using 
Theorem~\ref{universalgraphlemma}
along with a few standard techniques which we will not describe here, we can
find subsets $V_{i_1}, \ldots , V_{i_t}$ (for an appropriately chosen constant 
$t$), such that $d(V_{i_p},V_{i_q})>.999$
for every $1 \leq p <q \leq t$. This means that for  every 
$1 \leq p <q \leq t$, vertices of $V_{i_p}$ tend to be adjacent to a 
substantial proportion of the vertices of $V_{i_q}$.
On the other hand, if a substantial subset of $V_{i_p}$ is complete to 
a substantial subset of $V_{i_q}$, then we can apply induction to get a
large transitive subtournament in $T$, and so we may assume that no such 
subsets exist. We now construct a copy of $G$ in $T$, choosing at most one 
vertex from each of $V_{i_1}, \ldots, V_{i_t}$, and using the fact that
for $ 1 \leq p < q \leq t$ no substantial subset of $V_{i_p}$ is complete to a 
substantial subset of $V_{i_q}$ to obtain the backward edges in the galaxy 
ordering of $G$, thus obtaining the result of Theorem~\ref{galaxy}.

Obviously, every tournament obtained from a transitive tournament by adding a 
vertex is a galaxy.  It is not difficult to check that there is only one
tournament on five vertices that is not a galaxy. Here it is: its
vertex set is $\{v_1, \ldots, v_5\}$, and $v_iv_j$ is an edge if and only if
$(j-i) \; mod \; 5 \in \{1,2\}$. We call this tournament $S_5$.
We remark that $S_5$ is an example of a tournament that is obtained from a 
transitive tournament by adding two vertices, and that is not a galaxy.

Another result of \cite{BCC} is that:
\begin{theorem}
\label{S_5}
The tournament $S_5$ has the Erd\"{o}s-Hajnal property.
\end{theorem}
The proof of Theorem~\ref{S_5} uses similar ideas to the ones in the proof 
of Theorem~\ref{galaxy}.   Theorem~\ref{galaxy} and Theorem~\ref{S_5} together 
imply:
\begin{theorem}
\label{smalltourn}
Every tournament on at most five vertices has the Erd\"{o}s-Hajnal property.
\end{theorem}

We finish this section with another  curious corollary of Theorem~\ref{galaxy}.
Let $P_{k}$ denote a tournament of order $k$ whose vertices can be ordered 
so that the graph of backward edges is a $k$-vertex path.
\begin{theorem}
\label{path}
For every $k$, the tournament $P_{k}$ has the Erd\"{o}s-Hajnal conjecture.
\end{theorem}
Theorem~\ref{path} follows from the fact that, somewhat surprisingly,  $P_k$ has
a galaxy ordering. 

\section{Linear-size cliques, stable sets and transitive subtournaments}

In Conjecture~\ref{EHconj} and Conjecture~\ref{EHtourn}, every graph (or tournament) is conjectured 
to have a certain constant, lying in the $(0,1]$ interval, associated with it. 
A natural question is: when is this constant at its extreme? Excluding which 
graphs (or tournaments) guarantees a linear-size clique or stable set (or 
transitive subtournament)?

For undirected graphs this question turns out not to be interesting, because
if for some  graph $H$ every $H$-free  graph were to  contain a linear size 
clique or stable set, then $H$ would need to have at most two vertices
(we explain this later).
However, for tournaments the answer is quite pretty. We say that a tournament
$S$ is a {\em celebrity} if there exists a constant $0<c(S)\leq 1$ such that 
every
$S$-free tournament $T$ contains a transitive subtournament on at least 
$c(S)|V(T)|$ vertices. So the question is to  describe all celebrities. This
was done in \cite{heroes}, but before stating the result we need some 
definitions.

For a tournament $T$ and $X \subseteq V(T)$, we denote by $T|X$ the 
subtournament of $T$ induced by $X$.
Let $T_k$ denote the transitive tournament on $k$ vertices. If $T$ is a tournament and $X,Y$ are disjoint subsets of $V(T)$ such that $X$ is complete to $Y$,
we write $X\Rightarrow Y$. We write $v\Rightarrow Y$ for $\{v\}\Rightarrow Y$, and $X\Rightarrow v$ for $X \Rightarrow \{v\}$.
If $T$ is a tournament and $(X,Y,Z)$ is a partition of $V(T)$ into nonempty 
sets satisfying $X\Rightarrow Y$, $Y\Rightarrow Z$, and $Z\Rightarrow X$, we 
call $(X,Y,Z)$ a {\em trisection} of $T$. If $A,B,C,T$ are tournaments, and 
there is a trisection $(X,Y,Z)$ of $T$ such that 
$T|X,T|Y,T|Z$ are isomorphic to $A,B,C$ respectively, we write 
$T = \Delta(A,B,C)$. A {\em strongly connected component} of a tournament is
a maximal subtournament that is strongly connected.
One of the main results of \cite{heroes} is the
following:

\begin{theorem}
\label{celebrity}
A tournament is a celebrity   if and only if all its strongly connected 
components are celebrities. A strongly connected  tournament with more than 
one vertex is a celebrity if and only if
it equals $\Delta(S,T_k,T_1)$ or $\Delta(S, T_1, T_k)$ for some celebrity $S$ 
and some integer $k\ge 1$.
\end{theorem}

Following the analogy between stable
sets in graphs and transitive subtournaments in tournaments, let us 
define the {\em chromatic number} of a tournament $T$ to be the smallest
integer $k$, for which $V(T)$ can be  covered by  $k$ transitive 
subtournaments of $T$. We denote the chromatic number of $T$ by $\chi(T)$.
Here is a related concept: let us say that a tournament $S$ is a 
{\em hero} if there exists a constant $d(S)>0$ such that every $S$-free 
tournament $T$ satisfies $\chi(T) \leq d(S)$. Clearly every hero is a 
celebrity. Moreover, the following turns out to be true (see \cite{heroes})

\begin{theorem}
\label{hero}
A tournament is a hero if and only if it is a celebrity.
\end{theorem}

Another result of~\cite{heroes} is a complete list of all minimal non-heroes
(there are five such tournaments).

Let us now get back to undirected graphs. What if instead of asking
for excluding a single graph $H$ to guarantee a linear size clique or stable 
set, we ask the same question for a family of graphs? 
Let us say that a family $\mathcal{H}$ of graphs is {\em celebrated} if
there exists a constant $0<c(\mathcal{H}) \leq 1$ such that
every $\mathcal{H}$-free graph $G$ contains either a clique or a stable set of
size at least $c(\mathcal{H})|V(G)|$.
The {\em cochromatic number} of a graph $G$ is the minimum number of stable sets and cliques with 
union $V(G)$. We denote the cochromatic number of $G$ by $co\chi(G)$.
Let us say that a
family $\mathcal{H}$ is {\em heroic} if there exists a constant 
$d(\mathcal{H})>0$ such that  $co\chi(G) < d(\mathcal{H})$  for every
every $\mathcal{H}$-free graph $G$. Clearly, if $\mathcal{H}$ is heroic, then
it is celebrated. Heroic families were studied in \cite{CSeymour}.

Let $G$  be a   complete  multipartite  graph with $m$ parts, each of size 
$m$. Then $G$  has $m^2$  vertices, and no clique or stable set of size larger 
than $m$;  and the  same is true for $G^c$. Thus
every celebrated family contains a complete multipartite graph and
the complement of one. Recall that for every positive integer $g$ there exist 
graphs  with girth at least $g$ and no linear-size stable set (this is a 
theorem of Erd\"{o}s \cite{Erdos2}).  Consequently, every
celebrated family must also contain a graph of girth at least $g$, and,
by taking complements, a graph whose complement has girth at least $g$.
Thus, for a finite family of graphs to be celebrated, it must contain a
forest and the complement of one. In particular, if a celebrated family
only contains one graph $H$, then $|V(H)| \leq 2$.
The following conjecture is proposed in 
\cite{CSeymour}, stating that these necessary conditions for a finite family 
of graphs to be  celebrated are in fact sufficient for being heroic.

\begin{conjecture}
\label{heroicconj}
A finite family of graphs is heroic if and only if it contains
a complete multipartite graph, the complement of a complete multipartite graph,
a forest, and the complement of a forest.
\end{conjecture}

We remark that this is an extension of a well-known conjecture
made independently by Gy\'{a}rf\'{a}s \cite{Gyarfas} and Sumner \cite{Sumner},
that can be restated as follows in the language of heroic families:

\begin{conjecture}
\label{Gyarfas}
For every complete graph $K$ and every tree $T$, the family $\{K,T\}$ is heroic.
\end{conjecture}

The main result of \cite{CSeymour} is that Conjecture~\ref{heroicconj} and 
Conjecture~\ref{Gyarfas}
are in fact equivalent. Since a complete graph is a multipartite graph, 
the complement of one, and the complement of a forest, we deduce that
Conjecture~\ref{heroicconj} implies Conjecture~ \ref{Gyarfas}. The converse is a consequence of the
following theorem of \cite{CSeymour}:

\begin{theorem}
\label{ksplit}
Let $K$ and $J$ be graphs, such that both $K$ and $J^c$ are complete 
multipartite. Then there exists a constant $c(K,J)$ such that 
for every $\{K,J\}$-free graph $G$, $V(G)=X \cup Y$, where
\begin{itemize}
\item $\omega(X) \leq c(K,J)$, and
\item $\alpha(Y)  \leq c(K,J)$.
\end{itemize}
\end{theorem}

The situation for infinite heroic families is more complicated. Another
open conjecture of Gy\'{a}rf\'{a}s \cite{Gyarfas2} can be restated to say that 
a certain infinite family of graphs is heroic:
\begin{conjecture}
\label{hole}
For every complete graph $K$, and every integer $t>0$, the family consisting
of $K$  and all cycles of length at least $t$ is heroic. 
\end{conjecture}

If Conjecture~\ref{hole} is true, this is an example of a heroic set that does not 
include a minimal heroic set.

\section{Near-linear transitive subtournaments}

In this section we discuss an extension of the property of being a hero
studied in \cite{pseudo}.
Let us say that $\epsilon \ge 0$ is an {\em EH-coefficient} for a tournament $S$ if there exists $c>0$ such that 
every $S$-free tournament $T$ satisfies $\alpha(T)\ge c|V(T)|^{\epsilon}$.
(We introduce $c$ in the definition of the Erd\"{o}s-Hajnal coefficient to 
eliminate
the effect of tournaments $T$ with bounded number of vertices; now, whether 
$\epsilon$ is an EH-coefficient for $S$ depends only on arbitrarily large 
$S$-free tournaments.) Thus, Conjecture~\ref{EHconj} is equivalent to:
\begin{conjecture}
\label{EHcoeff}
Every tournament has a positive EH-coefficient. 
\end{conjecture}
If $\epsilon$ is an EH-coefficient for $S$, then so is every smaller non-negative number; and thus a natural invariant is the supremum of the set of all EH-coefficients for $S$.
We call this the {\em EH-supremum} for $S$, and denote it by $\xi(S)$. We 
remark  that the EH-supremum for $S$ is {\em not} necessarily itself an 
EH-coefficient  for $S$ (we will see an example later). One of the results 
of~\cite{pseudo} is a characterization
of all tournaments with EH-supremum $1$; and not all of these tournaments
turn out to be celebrities (in this language, a celebrity is a tournament
for which $1$ is its EH-coefficient, and not just its EH-supremum). 

The following theorem from \cite{Choromanski} suggests that EH-suprema
tend to be quite small:

\begin{theorem}
\label{upper}
Let $\mathcal{H}^{n,c}$ be the set of all $n$-vertex tournaments having 
EH-supremum at most $\frac{c}{n}$, where $c$ is an arbitrary constant such 
that $c > 4$, and 
let $\mathcal{H}^n$be the set of all
$n$-vertex tournaments. Then 
$${\frac{|\mathcal{H}^{n,c}|}{|\mathcal{H}^{n}|}} \rightarrow 1$$
as $n \rightarrow \infty$.
\end{theorem}

We say that a tournament $S$ is 
\begin{itemize}
\item a {\em pseudo-hero} if there exist constants $c(S),d(S) \ge 0$ such that every 
$S$-free  tournament $T$ with $|V(T)|>1$ satisfies 
$\chi(T) \leq c(S) (\log(|V(T)|))^{d(S)}$; and 
\item a {\em pseudo-celebrity} if there exist constants $c(S) >0$ and 
$d(S)\ge 0$ such that every $S$-free tournament $T$ with $|V(T)|>1$ satisfies 
$\alpha(T)\ge  c(S)\frac{|V(T)|}{(\log(|V(T)|))^{d(S)}}$.
\end{itemize}

In~\cite{pseudo} all pseudo-celebrities and pseudo-heroes are described 
explicitly.

\begin{theorem}
\label{pseudo}
The following statements hold:
\begin{itemize}
\item A tournament is a pseudo-hero if and only if it is a pseudo-celebrity.
\item A tournament is a pseudo-hero if and only if all its strongly connected 
components are pseudo-heroes.
\item A strongly-connected tournament with more than one vertex is a pseudo-hero if and only if either
\begin{itemize}
\item it equals $\Delta(T_2,T_k,T_l)$ for some $k,l \geq 2$, or
\item it equals $\Delta(S,T_1,T_k)$ or $\Delta(S,T_k,T_1)$ for some pseudo-hero $S$ and some integer $k>0$.
\end{itemize}
\end{itemize}
\end{theorem}

We remind the reader that by Theorem~\ref{celebrity} the tournament  
$\Delta(T_2,T_2,T_2)$ is not a celebrity, and yet by Theorem~\ref{pseudo} it is 
a pseudo-celebrity. Thus it is an example of a tournament that does not
attain its EH-supremum as an EH-coefficient.

We conclude this section with another result from~\cite{pseudo}, that shows
that after $1$, there is a gap in the set of EH-suprema.

\begin{theorem}
\label{gap}
Every tournament $S$ with $\xi(S)>5/6$ is a pseudo-hero and hence satisfies 
$\xi(S) = 1$.
\end{theorem}

The reason for  Theorem~\ref{gap} is another theorem from~\cite{pseudo} that
states that a tournament is a pseudo-hero if and only if it is 
$\mathcal{S}$-free for a certain family $\mathcal{S}$ consisting of  
six tournaments ($S_5$ is one of them). Thus for any tournament $T$ that is not
a pseudo-hero, $\xi(T) \leq \max_{s \in \mathcal{S}} \xi(S)$, and it is
shown that $\max_{s \in \mathcal{S}} \xi(S) \leq \frac{5}{6}$.

\section{Acknowledgment}

The author would like to thank Nati Linial for introducing her to the
world of the Erd\"{o}s-Hajnal conjecture at an Oberwolfach meeting  a number of
years ago. We are also grateful to Irena Penev and Paul Seymour  for their 
careful reading of the  paper, and for many valuable suggestions.

\end{document}